%&amstex
%   Warning   ........
%   amsppt.sty is necessary on top of amstex.tex
%                 ^^^^^^^^^
%\input amstex.tex
\input amsppt.sty
% double magnification
% \magnification=\magstep2
% single magnification
\magnification=\magstep1
\def\Mag@#1#2{\ifdim#1<1pt\multiply#1 #2\relax\divide#1 1000 \else
  \ifdim#1<10pt\divide#1 10 \multiply#1 #2\relax\divide#1 100\else
  \divide#1 100 \multiply#1 #2\relax\divide#1 10 \fi\fi}
\def\scalelinespacing#1{\Mag@\baselineskip{#1}\Mag@\lineskip{#1}%
  \Mag@\lineskiplimit{#1}}

\scalelinespacing{\magstephalf}

% \pageheight{6.9 true in} % this line is for the nec printer

\define\tcup {\text{$\tsize\bigcup $}}
\define\dcup {\bigcup}
\define\tcap {\text{$\tsize\bigcap $}}
\define\dcap {\bigcap}
\define\leftpar {\par\smallpagebreak\flushpar}

\newcount\onemore
\onemore=1
\define\mytag{\tag{\the\onemore}\global\advance \onemore by 1}

\nopagenumbers
\font\smalley=cmr8
\font\biggey=cmr12
% scaled\magstep1
% {\obeylines {\obeyspaces  }}

% \topmatter
% \vskip 6 truecm
% \endtopmatter

\document
\leftpar
{}
\leftpar
{}
\leftpar
{}
\leftpar
{}
\leftpar
{}
\leftpar
{}
\leftpar
{\biggey
SPACES OF LIPSCHITZ FUNCTIONS ON BANACH SPACES
}
\vskip 1 truecm
\leftpar
CHARLES STEGALL \ \ {\smalley
Institut f\"ur Mathematik, Johannes Kepler Universit\"at,
A-4040 Linz, Austria
}
%  \leftpar
%   Spring, 1991; Autumn, 1992
\vskip 1 truecm

\leftpar
A remarkable theorem of R. C. James is the following:
suppose that $X$ is a Banach space and $C \subseteq X$
is a norm bounded, closed and convex set such that
every linear functional $x^* \in X^*$ attains its
supremum on $C$; then $C$ is a weakly compact set.
Actually, this result is significantly stronger than this statement;
indeed, the proof can be used to obtain other
surprising results. For example,
suppose that $X$ is a separable
Banach space and $S$  is a norm separable subset
of the unit ball of $X^*$
such that for each $x \in X$ there exists $x^* \in S$
such that $x^*(x) = \|x\|$ then  $X^*$
is itself norm separable (this is in {\bf [R]}).
If we call $S$ a support set, in this case, with respect
to the entire space $X$, one can ask questions about
the size and structure of a support set, a support set
not only with respect to $X$ itself but perhaps with
respect to some other subset of $X$@. We analyze one particular case
of this as well as give some applications.
We begin by restating an old result, some
consequences of which are, perhaps, not well known.

\proclaim {Main Theorem} Let $X$ be a Banach space and $K$ a weak*
compact subset of the dual $X^*$. The following are equivalent:
\leftpar
(i) if $Y$ is any separable Banach space and $U : Y \to X$
is any operator, then $U^*(K)$ is a norm separable subset of
$Y^*$;
\leftpar
(ii) if $R : X \to C(K)$ is the canonical operator then
$R$ maps bounded subsets of $X$ into equimeasurable subsets of
$C(K)$;
\leftpar
(iii) the operator $R$ factors through an Asplund space
(a Banach space whose dual has the Radon-Nikodym property);
\leftpar
(iv) for any regular Borel (with respect to the weak* topology)
measure $\mu $  on $K$ and any $\epsilon > 0$
there exists a norm compact subset $K_\epsilon $ of $K$
such that $|\mu| (K_\epsilon ) > |\mu|(K) - \epsilon $;
\leftpar
(v) given any subset $S$ of $K$ and any $\epsilon > 0$
there exists an $x \in X$ and a real $r$
so that $\{x^* \in S : x^*(x) > r\}$ is not empty
and has norm diameter less than $\epsilon $.
\endproclaim

\leftpar
In {\bf [S5]} it is proved that if $K$ satisfies one of the above
then so does the weak* closed convex hull of $K$@.
An immediate consequence of the equivalence of
(i) and (iv) above is

\proclaim {Theorem} Let $U : X \to Y$ be an operator
and $K \subseteq Y^*$ a weak* compact set. Suppose that
$T^*$ is both a weak* and a norm homeomorphism on $K$@.
Then $K$ satisfies any one of the above if and only if
$T^*(K)$ does.
\endproclaim

\leftpar
Actually, we may formulate the result above as only requiring that
$T^*$ is one to one on $K$ and $(T^*)^{-1} $ is
(norm to norm) in the first Baire class
(or, for that matter, in any Baire class).
If (i) is false all of the others are false {\bf [S5]}.
The other directions, except, perhaps, those involving (v),
are essentially results of Smulian,
Dunford, Pettis, Phillips and Grothendieck;
see {\bf [G]}, {\bf [S5]} and {\bf [S6]}
for more detail.
Condition (v) is implicit in the results of Smulian.
What is somewhat newer is
that these conditions are equivalent to:
\leftpar
(vi) there exists a (norm to norm) first Baire class function
$\lambda : X \to X^*$ so that $\lambda (x) \in K$ and
$\lambda (x)(x) = \sup _{y^* \in K} y^*(x)$ for all $x \in X$@.

\leftpar
That (vi) follows from (v) can be obtained from
{\bf [JR]} and  see {\bf [S3]} for a transparent proof.
The first use of selection theorems in this connection is in
{\bf [F]}. That (vi) implies (i) is almost explicitly stated in
{\bf [R]} and is based on results of R. C. James
(see {\bf [S1]} for a discussion of James' results).
We give a proof that a somewhat
different version of (vi) implies (i).

\leftpar
Let $X$ be a Banach space and $K$ a weak* compact subset
of the unit ball of $X^*$. Suppose that $S$ is  a subset of $X^*$ that is
Lindel\"of in the weak* topology (in particular,
$S$ is norm separable) and $K \setminus S \ne \emptyset $.
Fix $y^* \in K \setminus S $.
Observe that
$$S \subseteq \dcup _{\|x\| = 1} \{x^* : |(y^* - x^*)x| > 0\}.$$
We may choose a sequence $\{x_n\}$ with $\|x_n\| = 1$
so that
$$S \subseteq \dcup _n \{x^* : |(y^* - x^*)x_n| > 0\}.$$
Define
$$\phi (x^*) = 1 - \sum _n 2^{-n}|(y^* - x^*)x_n|. $$
Observe that $-1 \le \phi (x^*)\le 1$ and
$$\sup _{z^* \in K} \phi (z^*) \ge 1 = \phi (y^*) > \phi (x^*)$$
if $x^* \in S$.
Thus, $\phi $ does not attain its supremum on $K$ at any point of $S$@.
Moreover,
$$ \align
&| \phi (x^*) - \phi (z^*) | =
\left| \sum _n 2^{-n} |(y^* - x^*)x_n| - \sum _n 2^{-n} |(y^* - z^*)x_n| \right|
\\
&\le \sum _n 2^{-n} |(z^* - x^*)x_n| \le \|z^* - x^*\|
\endalign $$
which proves that the Lipschitz norm of $\phi $ is no more than
one. Also, $\phi $ is concave and is
obtained by applying lattice operations
to the sequence $\{x_n\}$ and the constant functions.
Suppose, in addition, that $X$ is a separable Banach space.
If, given any $\phi $ that is weak* continuous on
$K$ and Lipschitz in the norm, there exists
$x^* \in S \cap K$ so that $\phi (x^*) = \sup _K \phi $
then $S \cap K = K$@. This follows from the fact that
the unit ball of $X^*$ is a separable metric space in the
weak* topology. If we fix $z^* \in K$ and let $S = K \setminus \{z^*\}$
then there is an $\phi $ that attains its supremum only
at $z^*$.
In fact, it suffices to take only those $\phi $ which can be
generated as above from $X$ and the constants.

\leftpar
Let $X^*$ be the dual of the Banach space $X$.
Let $L(X)$ denote the
smallest lattice of functions
that are both weak* continuous on
the unit ball $B$ of $X^*$  and
Lipschitz with respect to the norm on $X^*$; the
norm on $L(X)$ is
given by
$$L(f)= \max
\bigl\{
\sup _{\|x^*\| \le 1} |f(x^*)|,\ \sup
\Sb x^* \ne y^* \\ \|x^* \| \le 1 \\ \|y^* \| \le 1 \endSb
\|x^* - y^* \|^{-1} |f(x^*) - f(y^*)|
\bigr\}. \mytag $$
Observe that $L(X)$ is
generated by $X$ and the constants.
Consider the canonical operators
$$X @>I>>  L(X) @>R>> C(B).$$
This means that each $x^*$
in the unit ball of $X^*$ can be considered as three
different objects: $\delta _{x^*}$, denoting the point mass
in $C(B)$,
as an element of $L(X)^*$ and $x^* \in X^*$.
Also, $I$ is an isometry on $X$ and
$I^*$ is an isometry on $B$@.
Of course, $B$, considered as a subset of $L(X)^*$, does not have
the same convexity structure as it does as a subset of $X^*$.
Suppose that $K \subseteq B$ is weak* compact and satisfies
the hypothesis of the Main Theorem.
The canonical operator $R: L(X) \to C(K)$
(the restriction operator) has norm no more than one
and transforms the unit ball
of $L(X)$ into an equimeasurable subset of $C(K)$; this is nothing
more difficult than the results of {\bf [S6]} and the Arzela-Ascoli theorem.
There exist (see {\bf [JR]}, {\bf [F]} and {\bf [S3]})
a sequence of norm to norm continuous functions
$\lambda _n : L(X) \to L(X)^*$ such that
$\|\lambda _n (f)\| \le 1$ for all $n$ and $f \in L(X)$,
$\lim _n \lambda _n (f) = \lambda _0 (f) $  (limit in norm) exists
for all $f$,
$\lambda _0 (f) \in K$ and $f(\lambda _0 (f)) = \sup _K f$@.
Since $B$ is a retract of $X^*$ and $I^*$ is an isometry on
$B$, considered as a subspace of $L(X)^*$, it follows
that $B$ is a retract of $L(X)^*$.
We may assume that each $\lambda _n$ takes its values in $B$@.
%mark1
Define
$$S = \overline
{\dcup _{n = 0}^{\infty } \lambda _n (L(X))}^{\text {\, norm}}.$$
What is critical here is that $K \subseteq S$@.
Since $K$ and $S$ are both subsets of $B$
it does not matter whether we write $x^*$ or $\delta _{x^*}$.
Suppose that there exist $z^* \in K$ and
and $\eta > 0$ so that $\|z^* - x^*\| > \eta > 0$ for
all $x^* \in S$@.
Fix $\delta > 0$ and a norm compact subset $D$ of $S$.
Choose a finite subset $F$ of the unit ball of $X$ such that
$$\eta < \|x^* - z^*\| < \sup _{x \in F} |(x^* - z^*)x| + \delta $$
for all $x^* \in D$@.
Define
$$\phi (u^*) =
\sup \limits_{x \in F}
|(z^* - u^*)x|.$$
Now, the elementary computation:
$$ \align
&\phi (u^*) =
\sup \limits_{x \in F} |(z^* - u^*)x|
\le \sup \limits_{x \in F} \{|(v^* - u^*)x| + |(z^* - v^*)x| \}\\
&\le \sup \limits_{x \in F} \{\|v^* - u^*\| + |(z^* - v^*)x|\}
\le \|v^* - u^*\| + \phi (v^*).
\endalign $$
Since $\sup_{x^* \in K} \phi (x^*) \le 2$,
this proves that $L(1 - \phi ) \le 1$. Also, we have
$$ \align
&1 - \phi (u^*) \le 1 - \eta + \delta
\qquad \text {for all} \qquad u^* \in D \qquad \text {and}  \\
&1 - \phi (z^*) = 1.
\endalign $$
Let $\delta = \eta /3$.
We have produced a filter (over the norm compact subsets $D$ of $S$)
$$\Cal F = \{\psi _D : \psi
_D(D) \le  \frac {2} {3}\eta , \ \psi _D (z^*) = 1 ,
\text { and } L(\psi _D) \le 1\}.$$
Moreover, each $\psi _D$ is concave and generated by the lattice
operations on a finite subset of $X$ and the constants. Let $h$
be any cluster point  of $\Cal F$ in the topology of pointwise
convergence on $K$. If $F \subseteq X$ then we denote by $L(F)$
the smallest norm closed, linear sublattice of $L(X)$ containing
$F$ and the constants; if $F$ is norm separable then so is
$L(F)$. Choose $\psi _1 \in \Cal F$ so that
$|h(z^*) - \psi _1(z^*)| < \frac {\delta }{2}$.
Choose $F_1$ a finite subset
of $X$ such that $\psi _1 \in L(F_1)$ and let
$$L(F_1) =
\overline {\dcup _j G_{1,j}} $$
where each $G_{1,j}$ is norm
compact. In general, we construct $\{\psi _k \} \subseteq \Cal
F$, $\{F_k\}$ finite subsets of $X$, compact subsets
$\{G_{k,j}\}$ of $L(\cup _{i \le k} F _i)$ so that $\psi _k (z^*)
= 1 $ and
$$ \align
&L(\cup _{i \le k} F _i)) = \overline
{\dcup _j G_{k,j}} \quad \text { and } \\
&\psi _{k+1} \le \frac {2}{3} \eta  \quad \text { on }
\dcup _{p, q, r \le k} \lambda _p (G_{q,r}).
\endalign $$
Since
$\{\psi _k \}$ is bounded in the Lipschitz norm it follows that
the subset of $S$ where $\{\psi _k \}$ converges to null is norm
closed. Everything reduces to the case of the separable Banach
space $Y$ generated by $\cup _k F_k$. The functions $\{\psi _k
\}$ may be considered as weak* continuous and Lipschitz on the
unit ball of $Y^*$ and converge to null on every $\lambda _0(f)
| Y$ but this set contains all of $K | Y$ and we have that $\limsup
\psi _n (x^*) \le  \frac {2}{3} \eta  $
for all $x^* \in K$ which contradicts
$\liminf \psi _n (z^*) = 1$. Denote by $\Cal C \subseteq L(X)$
the set of concave functions such that $L(f) \le 1$.

\proclaim{Theorem} If $I : X \to L(X)$ is as above
and $K$ is a Radon-Nikodym set in the unit ball of $X^*$ then
$K = \overline { \lambda _0 (\Cal C) }^{\,\text {norm}}$.
\endproclaim

\demo {Proof} Suppose not. Suppose that $z^* \in K$, $\epsilon > 0$
and $\|\lambda _0 (f) - z^*\| > \epsilon $ for all $f \in \Cal C$@.
Suppose that $\|z^*\| < 1$. Choose $0 < \delta < \epsilon $ such that
$\|z^*\| + \delta < 1$. Then $\lambda  _0$ is a function of the
first Baire class with values in $C = B \setminus B(z^*, \delta )$.
The relevant geometric property of $C$ is that it is arcwise
connected and locally arcwise connected; see  {\bf [Ve]} and
its references or it is easy to rearrange the proof given in
{\bf [S3]} to obtain continuous functions $\xi _n : \Cal C \to C$
such that $\lim  \xi _n (f) = \lambda _0 (f)$ for all
$f \in \Cal C$@. We have already proved that
$$K \subseteq \overline {\dcup _n \xi _n (\Cal C)}^{\, \text {norm}}$$
but this contradicts $\|\xi _n (f) - z^*\| \ge \delta $
for all $n$ and all $f \in \Cal C$@. Suppose that $\|z^*\| = 1$.
Define $\rho : B \to B$ by $\rho (x^*) = x^*$ for
$\|x^* - z^*\| \ge \frac {\epsilon } {2} $ and
$$
\rho (x^*) = (1 -  \frac {\epsilon } {2} + \|x^* - z^*\|)x^*
+ (\frac {\epsilon } {2} - \|x^* - z^*\|)(-z^*)
$$
for $\|x^* - z^*\| \le \frac {\epsilon } {2}$.
Obviously, if $\|x^* - z^*\| \ge \frac {\epsilon } {2}$
then $\|\rho (x^*) - z^*\| \ge \frac {\epsilon } {2}$;
suppose that
$$\delta =  \frac {\epsilon } {2} - \|x^* - z^*\|
< \frac {\epsilon } {2}.$$
If
$\delta \le  \frac {\epsilon } {8} $ we have that
$ \|x^* - z^*\| \ge \frac {\epsilon } {2} - \frac {\epsilon } {8}$
and
$$\|\rho (x^*) - z^*\| = \|(1 - \delta )x^* - (1 + \delta )z^*\|
\ge  \|x^* - z^*\| - \delta \| x^*  + z^* \| \ge
  \frac {\epsilon } {2} - \frac {\epsilon } {8} -  \frac {2 \epsilon } {8}
=  \frac {\epsilon } {8}.
$$
If  $\delta \ge \frac {\epsilon } {8}$ then
$$
\|\rho (x^*) - z^*\| = \|(1 - \delta )x^* - (1 + \delta )z^*\|
\ge  (1 + \delta ) -  (1 - \delta ) = 2\delta \ge   \frac {2\epsilon } {8}.
$$
In either case,  $ \|\rho (x^*) - z^*\| \ge   \frac {\epsilon } {8}$.
Observe that $\lim _n \rho \circ \lambda _n (f) = \lambda _0 (f)$.
Again, we know that
$$K \subseteq
\overline {\dcup _n  \rho \circ \lambda _n  (\Cal C)}^{\, \text {norm}}$$
but this is a contradiction. In the formula for $\rho$ we have considered
$B$ as a subset of $X^*$; of course, when $B$ is considered as a subspace
of $L(X)$ then we have to consider the canonical inverse of the
isometry composed with $\rho$.
\enddemo

\leftpar
With a non trivial expansion of these techniques, the following
can be shown: with the hypothesis as above, if $x^* \in K$
and $\epsilon > 0$ then there exits $\phi \in \Cal C$ such that
$\phi $ attains its supremum only on $K \cap B(x^*,\epsilon )$,
which means that $\lambda _0 (\phi ) \in  K \cap B(x^*,\epsilon ) $
(for any $\lambda _0$).
The case considered above where $\|z^*\| < 1$ could be dealt
with differently, with considerable overkill. In any infinite
dimensional Banach space the sphere $S$ of a closed ball $B$
is a retract of $B$; this is an old result for a Hilbert space and
it follows for any Banach space by the Anderson-Kadec-Torunczyk
theorem (see {\bf [BP]} and {\bf [vM]}). Thus we could retract
$B$ onto $B \setminus B(z^*, \delta )$ and replace
each $\lambda _n$ with its composition with such a retraction.
The following is almost proved in {\bf [F]}; what is missing
can be found in {\bf [R]}.
More details can be found in
{\bf [FG]} and short proofs in {\bf [S9]} and {\bf [S10]}.
The proof below does not use the results of R. C. James.

\proclaim {Theorem} Let $X$ be a Banach space such that
$X^*$ has the Radon-Nikodym property (the unit ball of
$X^*$ satisfies the Main Theorem).
Then there exists an interval of ordinals $[1, \eta ]$
such that for each $1 \le \alpha \le \eta  $ we may construct
subspaces
$X_\alpha $ of $X$ and subspaces $X_\alpha ^\sharp$
of $X^*$ such that
\leftpar
(i) $X_\alpha ^\sharp | X_\alpha $ is an isometry onto $X_\alpha ^*$;
\leftpar
(ii) for each limit ordinal $\beta $ it follows that
$$
X_\beta  = \overline { \dcup X_\alpha }^{\,\text {norm}}
\qquad \text {and} \qquad
X_\beta ^\sharp = \overline { \dcup X_\alpha ^\sharp }^{\,\text {norm}};
$$
\leftpar
(iii) $X_1$ is norm separable;
\leftpar
(iv) the norm density of $X_\alpha $
(which equals the norm density of $X_\alpha ^\sharp $)  is countable or
no more than that of $\alpha $ and
\leftpar
(v) $X_\eta = X$.
\leftpar
\endproclaim

\demo {Proof} Here, $K$ is the unit ball of $X^*$.
The proof is by induction. Assume that we have constructed
$Y_1 = X_\alpha $.
Let $L_i$  be the smallest lattice in $L(X)$  containing $Y_i$.
Choose $Y_1 \subseteq Y_2 $ such
that $Y_2 \ne Y_1 $,
$$X_\alpha ^\sharp \subseteq \big[ \dcup _j \lambda _j(L_2) \big]
\subseteq X^*$$
and both $Y_2$ and $Y_1$ have the same norm density.
In general, choose $Y_n \subseteq Y_{n+1}$ so that
$Y_{n+1}$ norms $[\cup _j \lambda _j(L_n)]$.
Remember that the mapping from $L(X)^*$ to $X^*$ is an isometry on the
unit ball of $X^*$. Define
$X_{\alpha +1} = [\cup _n Y_n]$ and
$X_{\alpha +1}^\sharp = [\cup _n \cup _j \lambda _j (L_n)]$.
It follows that $X_{\alpha +1}^\sharp | X_{\alpha +1} $
is an isometry onto.
\enddemo

\leftpar
Before returning to Banach spaces, we  discuss
an abstraction of the above.
Suppose that $K$ is a compact Hausdorff space that has a finer
topology given by a metric $d$@. We denote by
$L(K,d)$ the Banach space of functions that are both continuous on $K$
and Lipschitz with respect to $d$ and the norm on $L(K,d)$ is
given by
$$L(f) = \max
\bigl\{ \sup _k |f(k)|,\ \sup _{k \ne k^{\prime }} d(k,k^{\prime })^{-1}
|f(k) - f(k^{\prime }| \bigr\}. \mytag $$
We make two assumptions:
\leftpar
(i) $L(K,d)$ separates the points of $K$ and
\leftpar
(ii) for each regular Borel measure $\mu $ on $K$ and each $\epsilon > 0$
there exists a subset $M$ of $K$ that is compact in the metric and
$$|\mu |(M) > \|\mu \| - \epsilon .$$
\leftpar
If we replace $d$ by the metric
$$\rho (k, k^{\prime}) = \sup _{L(f) \le 1} |f(k) - f(k^{\prime })| \mytag $$
it is quite easy to see that $L(K,\rho) $ also satisfies (i) and (ii).
It may happen of course that $L(K,\rho) $  is larger than $L(K,d)$
and this, in general, is desirable.
Thus, we shall assume that $\rho = d$.
If we assume that $(K,d)$ is complete and
\leftpar
(iia) for each subset $S$ of $K$
and each $\epsilon > 0$ there exists an open subset $V$
(in the topology that makes $K$ compact) so that $V \cap S \ne \emptyset $
and diameter $V \cap S < \epsilon $
\leftpar
then (iia) implies (ii) ({\bf [S7]}) and
(i) and (iia) together give
a useful topological
characterization of compact spaces that are homeomorphic to
a weak* compact subspace of a dual Banach space satisfying
the Main Theorem (see  {\bf [S7]} and {\bf [Na]}). For our purposes,
and probably all
purposes, (i) and (ii) are more useful than
(i) and (iia). Condition (iia) does not imply condition (i).
Indeed, the examples in {\bf [T]}  and {\bf [SL]}  are compact
spaces satisfying (iia)
(this is a special case of {\bf [Gu]}, {\bf [V]}, {\bf [Na]}
and for lots of details, {\bf [S7]}) and this example
is not an Eberlein compact but is a Corson compact
(same references). As we shall see below, a Corson
compact space that is the continuous image of a
space satisfying (i) and (ii) must be an Eberlein
compact. Conditions (i) and (ii) together are much
stronger than (iia) alone. Observe that the canonical
operator $U : L(K,\rho) \to C(K)$ transforms the unit ball
of $L(K,\rho)$ into an equimeasurable set (see {\bf [S6]})
and it follows that $K$ considered as a subset of $L(K,\rho)^*$
satisfies the Main Theorem.
Suppose that there exists a subset $S$ of $L(K,d)$ such that $S$ also
separates the points of $K$ and $S$ is point countable on $K$,
which means that
$$\{f \in L : |f(k)| > 0\}$$
is countable for each $k \in K$ (this is the definition that
$K$ is a Corson compact).
We may also assume that
$$\rho (k, k^{\prime}) = \sup \Sb L(f) \le 1 \\ f \in S \endSb
|f(k) - f(k^{\prime })| \mytag $$
We shall show that if a compact Hausdorff space $K$ is both a Corson
compact space and homeomorphic (in the weak* topology) to a weak*
compact subset of some $X^*$ satisfying the conditions of the
Main Theorem then $K$ is an Eberlein compact.

\leftpar
Suppose that $X$ is a Banach space and $K \subseteq B \subseteq X^*$
are as above and $K$ satisfies the Main Theorem. Suppose that
$K$ is also a Corson compact. Define $L(X)$ as above and let
$$X @> I >> L(X) @> R >> C(K)$$
be the canonical operators; $R$ is the restriction operator which
is also a lattice homomorphism. Suppose that $S \subseteq C(K)$
is point countable and point separating. It is easy to see that the
smallest algebra over the rational numbers is also
point countable and point separating. Thus we may assume that
$S$ is uniformly dense in $C(K)$; also $R(L(X))$ is uniformly
dense in  $C(K)$ (see {\bf [N]} concerning the Stone-Weierstra\ss \
theorem). Suppose $S = \{f_\alpha \}$ and for each $\alpha $
choose an $n \in \Cal N$ and $g _{n, \alpha }$
in the unit ball of $L(X)$ such that
$$\| 2^n R(g _{n, \alpha }) - f_ \alpha \| < \frac {1} {4}$$
and define
$$h _{n, \alpha } = (g _{n, \alpha }  \vee 2^{-n -1}) -  2^{-n -1}.$$
If, for some $k \in K$,  $h _{n, \alpha }(k) > 0$, then
$R(2^n g _{n, \alpha } )(k) > \frac {1}{2}$ and it follows that
$f_ \alpha (k) \ne 0$; thus, $\{ h _{n, \alpha }\}$ is point countable.
Suppose, $k$ and $k^\prime$ are distinct points of $K$. Choose
$f \in C(K)$ such that $\|f \| = 1$ and $f(k) = 1$ and
$f (k^\prime) = -1$. Choose $\alpha $  such that
$\|f _\alpha - f\| < \frac {1} {4} $. For the appropriate
$2^n g _{n, \alpha }$ we have that
$$ \| 2^n R(g _{n, \alpha }) - f \| < \frac {1}{2} $$
and it follows that $ 2^n R(g _{n, \alpha })(k) >  \frac {1}{2}  $
and  $ 2^n R(g _{n, \alpha })(k^\prime) < - \frac {1}{2}$.
Thus, $ h _{n, \alpha }(k) > 0$ and  $ h _{n, \alpha }(k^\prime) = 0$.
We have that $\{ h _{n, \alpha }\}$ is point countable and
point separating. We may assume the following: we have a Banach space
$Y$ (in this case,  $Y = L(X)$) and a Radon-Nikodym subset
$K$ of $Y^*$ and a subset $T$ ($T =  \{ h _{n, \alpha }\} $)
of $Y$ that is point countable and point separating on $K$@.
We shall apply the entire procedure to $L(Y)$; this
roughly corresponds to the double interpolation in {\bf [S10]}.
At this point, however, we shall simplify the
notation.

\leftpar
Suppose that $X$ is a Banach space,  $K \subseteq X^*$ a
Radon-Nikodym set, and $T \subseteq X$  is point countable
and point separating on $K$@.
We may also assume that the linear span $[T]$ of $T$ is $X$;
in fact we may assume that $T$ is a norm dense linear (over
the rationals) subspace of $X$@. The weak* closed convex hull
of $K$ is also a  Radon-Nikodym set (see {\bf [S6]}
but not the review thereof in \text {\it Mathematical Reviews\/}).
We check that the norm closed convex hull of $K$
is the same as the  weak* closed convex hull. Let $\mu $ be a Radon
probability measure on $K$;
the resolvent $r(\mu ) \in X^*$ is the linear functional
given by
$$r(\mu ) (x) = \int x\,d\mu. $$
The set of all resolvents is exactly the
weak* closed convex hull of $K$@.
For any $\epsilon > 0$
there exists a norm compact set $F \subseteq K$ such that
$\mu (F) > 1 - \epsilon $. The norm closed convex hull of $F$
is norm compact and contains $r(\nu )$ where
$ \nu  = (\mu (F))^{-1} \chi _F \mu$. For any $x \in X$
we have that
$$\align
&|r(\mu ) (x) - r (\nu) (x)| = |\int x \,d\mu - \int x\nu | \\
&\le \left( \int _{K \setminus F} |x|\,d\mu \right)
+ |1 -  \mu (F) ^{-1}| \left( \int _{F} |x|\,d\mu \right)
\le \big( \epsilon + |1 - \frac {1} {1 - \epsilon }| \big)\|x\|.
\endalign $$
It remains to verify that $T$ is point countable  on the
weak* closed convex hull of $K$@. Suppose that $y^*$ is in this hull;
then there exists a countable set $C \subseteq K$ such that
$y^*$  is in the norm convex hull of $C$@. If some $x \in T$
is nonzero on  $y^*$ then $x$ is nonzero on some point in $C$;
since $C$ is countable this proves that
$$\{x \in T : |y^* (x)| > 0\}$$
is countable. Finally, we may assume that $K$ is convex.
Define $I : X \to L(X)$ and ${\lambda _0, \lambda _1, \lambda _2, \dots}$
as before. Since $K$ is convex and norm closed it is a (norm) retract of
$X^*$ (see  {\bf [BP]}); also, $I^*$ is an isometry on $K$
(considered as a subset of $L(X)$) and it is easy to see that
$K$ is a retract  $L(X)^*$. We may assume that each $\lambda _n $ takes
values in $K$. An elementary property of Corson compacta is that given
an infinite set $T^\prime \subseteq T$ then there exist
$S \subseteq T$,  $T^\prime \subseteq S$ and a continuous (in this
case weak* continuous) retraction $r$ on $K$ such that
$f \circ r = f$ for all $f \in S$ and $f \circ r = 0$ for all
$f \in T \setminus S$; such sets $S$ are said to be good subsets
(see {\bf [Ne]}). The increasing union of good sets is also good.
The crucial facts about any
weak* compact
subset $F$ of $K$ and any family $S$ of continuous functions
that separate the points of $F$
are the following:
$$\text {norm density } F \le \text { weight } F
\le \text { cardinality } S$$
(the weight is  with respect to the weak* topology)
because $K$ satisfies
the hypothesis of the Main Theorem.
There exists an interval of ordinals $[0,\eta ]$ for which
we can construct good increasing subsets $S_\alpha $ of $T$
with the following properties:
\leftpar
(i) $S_0 = \emptyset $;
\leftpar
(ii) $S_1$ is countable;
\leftpar
(iii) $S_\alpha $ is countable or
the cardinality of $S_\alpha $ is
the cardinality of $\alpha$;
\leftpar
(iv) if $\beta $ is a limit ordinal then
$S_\beta = \cup _{\alpha < \beta } S_\alpha $;
\leftpar
(v) $S_\eta  = T$;
\leftpar
(vi) there exist continuous retractions $\{r_\alpha \}$
on $K$ such that $fr_\alpha = f$ if
$f \in S_\alpha $ and if $f \in T \setminus S_\alpha $ then
$fr_\alpha = 0$ and
\leftpar
(vii) if $L_\alpha $ denotes the smallest lattice containing
$S_\alpha $ and the constants then, for any limit ordinal $\beta $,
$$
\overline { \dcup _n \lambda _n (L_\beta  ) }^{\,\text {norm} }
= \overline
{\dcup _n \dcup _{\alpha <\beta } \lambda _n (L_\alpha ) }^{\,\text{norm}}
= r_\beta (K).
$$
The construction is by induction and follows from the following
discussion.
Fix an infinite subset $S$ of $T$. Construct
an increasing subsequence of
good subsets $\{T_j\}$
of $T$ of the same
cardinality as $S$, with associated retractions $\{r_j\}$,
the norm closure of $T_j$ is a linear subspace $X_j$ of $X$,
$X_j$  generates the lattice   $L_j$ and
$$r_j(K) \subseteq
\overline { \dcup _n \lambda _n (L_j) }^{\,\text {norm}}
\subseteq r_{j+1}(K).$$
It follows that $\cup _j T_j = T_0$ is also good with retraction $r$@.
Let $L = \overline {\cup _j L_j} $; it follows that
$L$ is the lattice generated by
$Y = \overline {\cup _j X_j} ^{\,\text {norm}}$.
Observe that
$f\circ r = f$ for all $f \in L$@.
Let $J : L \to L(X) $  be the containment operator. From the above we know
that
$$\align
&J^* (K) = J^* r(K)
\subseteq \overline {\dcup _n J^*\lambda _n (L) }^{\,\text {norm}} \\
&\subseteq \overline
{\dcup _j \dcup _n J^*\lambda _n (L_j) }^{\,\text {norm}}
\subseteq
\overline { \dcup _j J^*r_{j+1}(K) }^{\,\text {norm}} \subseteq J^*r(K).
\endalign $$
For a given $\epsilon > 0$ and $k \in K$ choose $n$, $j$
and $\phi \in L_j$ such that
$$\|J^* r(k) - J^* (\lambda _n (\phi )) \| < \epsilon ;$$
because $f \circ r \in L$ and
$r (\lambda _n (\phi )) =  \lambda _n (\phi )$ it follows that
$$\align
&\|r(k) -  \lambda _n (\phi )\|
= \sup \Sb f \in L(X) \\ L(f) \le 1 \endSb
|f(r(k)) -  f(\lambda _n (\phi )) | \\
&= \sup \Sb f \in L \\ L(f) \le 1 \endSb
|f(r(k)) -  f(\lambda _n (\phi )) |
 = \| J^* r(k) - J^* (\lambda _n (\phi )) \| < \epsilon .
\endalign
$$
The relevant formula is
$$r(K) =
\overline { \dcup _j r_j (K) }^{\,\text {norm}}.$$
If we have constructed $S_\alpha $ then $S_{\alpha +1}$ is constructed by
this process.

\leftpar
If we define
$$\rho _\alpha (k, k^{\prime}) =
\sup \Sb L(f) \le 1 \\ f \in S_\alpha \endSb
|f(k) - f(k^{\prime })| \mytag $$
then we have that
$$\rho _\alpha (k, k^{\prime}) =
\|r _\alpha (k) -  r _\alpha (k^{\prime}) \|.$$
Define the sets $T_\beta =
S_\beta \setminus \cup _{\alpha < \beta} S_\alpha $ and the functions
$$\phi _\beta  (k) =
\sup \Sb L(f) \le 1 \\ f \in T _\beta \endSb |f(k)|.$$
Observe that each $\phi _\beta$ is lower semi-continuous,
$\phi _\beta  (k) = 0$ if $\beta $ is a limit ordinal and
$$\phi _{\beta + 1} (k)
= \rho _{\beta + 1} (r_{\beta + 1}(k), r_ \beta (k))
=\|r_{\beta + 1} (k) -  r_ \beta (k) \|.$$
Suppose that each $\alpha _i$ is not a limit ordinal and
$$0 < \alpha _1 <  \alpha _2 < \dots < \beta = \sup _n \alpha _n < \eta .$$
Fix $k \in K$  so that $r_\beta (k) = k$.
Choose $k_n \in r_{\alpha _n} (K)$ so that $\|k - k_n \| \to 0$.
Then
$$ \align
&\|k_n - r_{\alpha _{n+1}-1}(k) \|
=\|r_{\alpha _{n+1}-1} (k_n) - r_{\alpha _{n+1}-1}(k) \| \\
&= \rho _{\alpha _{n+1}-1} (k_n,k) \le \|k_n -  k\| @>n>> 0.
\endalign $$
Therefore, $\|k - r_{\alpha _{n}-1}(k) \| @>n>> 0$ and
$$\phi _{\alpha _n } (k) \le
\| r_{\alpha _{n}-1}(k) - r_{\alpha _n}(k) \|
= \rho _{\alpha _n - 1}(k, r_{\alpha _n}(k)) @>n>> 0.$$
This proves that
$$\{\alpha : |\phi _\alpha (k)| > \epsilon \}$$
is finite for each $k$ and each $\epsilon > 0$.
Assume that for any subset $T$ of $S$ such that
the cardinality of $T$ is less than that of $S$
we have that the image of $K$ in $\prod _T [-1,1]$
defined by $k \to (f(k))_{f \in T}$ is an Eberlein compact.
Let
$$\bigl\{ E_\beta \subseteq \prod _{T _\beta }[-1,1] \bigr\}$$
be the images under these maps.
>From this, it follows that $K$ itself is an Eberlein
compact because the image of $K$ in $\prod E_\beta$ defined by
$k \to ((f(k))_{f \in T_\beta })_\beta $
is a homeomorphism whose image is in the $c_0$ sum
of the $E_\beta$'s; this follows because
$\{\alpha : |\phi _\alpha (k)| > \epsilon \}$
is finite for each $k$ and each $\epsilon > 0$.
The following is in {\bf [S10]} but
the proof here is even easier.

\proclaim {Theorem}
Suppose that
$K$ and $T$ are compact Hausdorff spaces and
$q : K \to T$ is continuous.
Suppose that at least one of $K$ and $T$ is
a Corson compact and at least one is a  Radon-Nikodym compact.
Then $q$ continuously factors through a space $S$ that is an
Eberlein compact.
\endproclaim

\demo {Proof} We have shown that if $K$, or $T$, has both properties
then $K$, or $T$, is an Eberlein compact. It is known that if
$K$ is Corson compact then $q(K)$ is Corson compact; see  {\bf [Ne]},
but a nearly trivial proof is given in  {\bf [S10]}, which bears
repeating. We shall assume that $T$ is Corson compact and
$K$ is a  Radon-Nikodym compact.
We assume that $q$ is onto and $A \subseteq C(K)$
is the subalgebra such that $f \in A$ if and only if
there exists $h \in C(T)$ such that $gq = f$.
We have already used a variant of the following.
Suppose that $C$ is a subset of the unit ball of
$C(K)$ that is equimeasurable and separates the points of $K$@.
We may assume that $C$ is a convex and symmetrical
subset of the unit ball with $1_K \in C$@.
It is easy to check that $C\cdot C$ is also equimeasurable
and, by induction, $C^n$ is equimeasurable;
it is also routine to check that $(C\vee C) - C$
is equimeasurable.
Let $E$ be an equimeasurable set that
is closed, convex, symmetrical and
contains
$$\sum _n 2^{-n} C^n.$$
Again, it follows from  the Stone-Weierstra\ss \ theorems
(see  {\bf [N]}) that $\cup _n 2^n E$ is uniformly
dense in $C(K)$.
Let ${\bold F}$ be any subset
of $A$ that is point countable
and separates the (positive and multiplicative)
states of $A$. We shall show that there exists a subset
${\bold G}$ of $C(K)$ that is both equimeasurable and point
countable and the algebra $A_1$ generated by ${\bold G}$
contains ${\bold F}$@. Observe that the family of all
polynomials in ${\bold F}$ with rational coefficients
is also point countable.
Thus, we may assume that
${\bold F}$ is a dense subset of the unit ball
of $A$ (yes, again, the Stone-Weierstra\ss \ theorem).
Partition ${\bold F}$ by
$${\bold F} = \dcup _n {\bold F}_n$$
so that
$${\bold F}_n \subseteq (2^n  E) + B(0, \frac {1}{4}).$$
For each $f_{n, \alpha } \in {\bold F}_n$
choose $h_{n, \alpha } \in E$ so that
$\|2^nh_{n, \alpha } - f_{n, \alpha }\| < \frac {1}{4}$.
Define
$$u_{n, \alpha } = (2^nh_{n, \alpha } \vee \frac {1}{2}) -
\frac {1}{2} \ge 0$$
and observe that
$$2^{-n}u_{n, \alpha } = (h_{n, \alpha } \vee \frac {1}{2^{n+1}})
- \frac {1}{2^{n+1}} \in (E \vee \frac {1}{2^{n+1}}) -
\frac {1}{2^{n+1}} \subseteq (E \vee E) - E$$
because $E$ is convex and contains the origin and $1_K$.
Thus, $\{2^{-n}u_{n, \alpha } : n, \alpha \}$ is equimeasurable.
Fix $k \in K$; if $2^{-n}u_{n, \alpha } (k) > 0$ then
$2^nh_{n, \alpha }(k) > \frac {1}{2}$ which implies that
$|f_{n, \alpha }(k)| > \frac {1}{4}$. Thus,
$\{2^{-n}u_{n, \alpha } : n, \alpha \}$ is also point countable.
Fix two points $k_0$ and $k_1$ in $K$ so that
$k_0|A = a_0$ and $k_1|A = a_1$ are distinct (positive,
multiplicative) states of $A$
and choose a function $h \in A$
so that
$$-1 = h(a_0) < h(a_1) = 1 = \|h\|.$$
There exist $n$ and $f_{n, \alpha } \in {\bold F}_n$ so that
$\|f_{n, \alpha } - h\| < \frac {1}{4}$. Therefore,
$\|2^nh_{n, \alpha } - h\| < \frac {1}{2} $
and it follows that
$$0 = u_{n, \alpha }(a_0) < u_{n, \alpha }(a_1).$$
Let ${\bold G} = \{2^{-n}u_{n, \alpha } : n, \alpha \}$.
Let $A_2$ be the algebra generated by
${\bold F} \cup {\bold G}$
and let $A_1$ be the algebra generated by
${\bold G}$.
The arguments above show that the
functions in
$A_2$ and $A_1$ separate exactly the same points of $K$@.
If some $h \in A$ separates two points in $K$
then those two points are also separated by some $u_{n,\alpha }$.
This means that
$A \subseteq A_1 = A_2$.
The algebra $A_1$ is generated by the point countable and equimeasurable
set $\{2^{-n}u_{n, \alpha }: n, \alpha  \}$.
We have proved that the state space
of $A_1$ is an Eberlein compact. From {\bf [Gu1]} or {\bf [BRW]}
(see also {\bf [MR]} for another approach)
it follows that
the state space of $A$ is also an Eberlein compact and this, of course, is
$q(K) = T$@.
\enddemo
%%%%%

% \heading
\leftpar
Bibliography
%\endheading
\leftpar
{\bf [AL]} D. Amir and J. Linden\-strauss, The structure of
weakly compact sets
in Banach spaces, Ann. of Math. 88 (1968), 35--46
\leftpar
{\bf [BRW]} Y. Benjamini, M. E. Rudin and M. Wage,
Continuous images of weakly compact subsets of Banach spaces,
Pacific J. Math. 70 (1977), 309--324
\leftpar
{\bf [BP]} C. Bessaga and A. Pelczynski, {\it Selected Topics in
infinite-dimensional Topology}, PWN Warszawa, 1975
\leftpar
{\bf [DFJP]} W. Davis, T. Figiel, W. Johnson and A. Pe\l czynski,
Factoring weakly compact operators, J. Funct.
Analysis 17 (1974), 311--327
\leftpar
{\bf [F]} M. Fabian, Each weakly countably determined Asplund
space admits a Fr\' echet differentiable norm, Bull. Austral.
Math. Soc. 36 (1987), 367--384
\leftpar
{\bf [FG]} M. Fabian and G. Godefroy, The dual of every Asplund space
admits a projectional resolution of identity, to appear in Studia Math.
\leftpar
{\bf [G]} A. Grothendieck, {\it Produits tensoriels
topologiques et espaces nucl\'eaires}, Memoire AMS 16,
1955
\leftpar
{\bf [Gu]} S. P. Gul'ko, On the structure of spaces of continuous
functions and their complete paracompactness,
Russian Math. Surveys 34 (1979), 36--44
\leftpar
{\bf [Gu1]} S.P. Gul'ko, On properties of subsets of $\sigma$-products,
Dokl. Akad. Nauk. USSR 237 (1977), 505--508
\leftpar
{\bf [Gr]} G. Gruenhage, A note on Gul'ko compact spaces, Proc. AMS
100, (1987) 371--376
\leftpar
{\bf [MR]} E. Michael and M. E. Rudin, A note on Eberlein compacts,
Pacific J. Math. 72 (1977), 487--495
\leftpar
{\bf [vM]} J. van Mill, {\it Infinite-Dimensional Topology},
North-Holland, Amsterdam, 1989
\leftpar
{\bf [N]} M. A. Naimark, {\it Normed Rings}, Noordhoff, Groningen,
1964
\leftpar
{\bf [Na]} I. Namioka, Radon-Nikodym compact spaces and fragmentability,
Mathematika 34 (1987), 258--281
\leftpar
{\bf [Ne]} Negrepontis, Banach Spaces and Topology, in {\it Handbook of Set
Theoretic Topology}, North Holland, Amsterdam, 1984
\leftpar
{\bf [R]} G. Rod\'e, Superkonvexit\"at und schwache Kompaktheit,
Arch. Math. 36 (1981), 62--72
\leftpar
{\bf [SL]} G. A. Sokolov and Leiderman, Adequate families of sets
and Corson compacts,
Commentationes Mathematicae  Universitatis Carolinae 25 (1984),
233--246
\leftpar
{\bf [S1]} C. Stegall, {\it Applications of Descriptive Topology in
Functional Analysis}, Universit\"at Linz, 1985
\leftpar
{\bf [S2]} C. Stegall, The duality between Asplund spaces and spaces with
the Radon-Nikodym property, Israel J. Math. 29 (1978),
408--412
\leftpar
{\bf [S3]} C. Stegall,
Functions of the first Baire class with values in Banach spaces,
Proc. AMS 111 (1991), 981--991
\leftpar
{\bf [S4]} C. Stegall,
A Proof of the Theorem of Amir and Lindenstrauss,
Israel J. Math. 68 (1989), 185--192
\leftpar
{\bf [S5]} C. Stegall,
The Radon-Nikodym property in conjugate Banach spaces,
Trans. AMS 206 (1975), 213-223
\leftpar
{\bf [S6]} C. Stegall, The Radon-Nikodym property in
conjugate Banach spaces II,
Trans. AMS 264 (1981), 507-519
\leftpar
{\bf [S7]} C. Stegall, The Topology of Certain Spaces of Measures,
Topology and its Applications 41 (1991), 73--112
\leftpar
{\bf [S8]} C. Stegall, Dunford-Pettis Sets, in {\it General Topology and
its Relations to Modern Analysis and Algebra VI\/},
Z. Frol\'\i k, editor, Heldermann, Berlin, 1988
\leftpar
{\bf [S9]} C. Stegall,
More Facts about conjugate Banach Spaces
with the Radon-Nikodym Property,
Acta Universitatis Carolinae--Math. et Phys. 31 (1990), 107--117
\leftpar
{\bf [S10]} C. Stegall,
More Facts about conjugate Banach Spaces
with the Radon-Nikodym Property II,
Acta Universitatis Carolinae--Math. et Phys. 32 (1991), 47--54
\leftpar
{\bf [T]} M. Talagrand, Espaces de Banach faiblement k-analytique,
Ann. of Math. 110 (1979), 407--438
\leftpar
{\bf [V]} L. Vasak, On one generalization of weakly compactly
generated Banach spaces, Studia Math. 70 (1981), 11--19
\leftpar
{\bf [Ve]} Libor Vesel\' y, Characterization of Baire-one
functions between topological spaces, to appear,
Acta Universitatis Carolinae--Math. et Phys.
\enddocument
\bye \bye